\documentclass{amsart}
\usepackage[a4paper,margin=3.2cm]{geometry}
\usepackage{graphicx} % Required for inserting images
\usepackage{amsmath,amsthm,tikz,tikz-cd,bbm,subcaption,url}
\usepackage{blkarray,booktabs,bigstrut}
\usepackage{mathtools}
\usepackage{verbatim}
\usepackage[url=false]{biblatex}
\usepackage[dvipsnames]{xcolor}
\title{Mutation Of Matrices Over Group Rings}
\author{Dani Kaufman}
\author{Carmen Alves Sabin}
\date{June 2025}

\usepackage{hyperref}
\hypersetup{
	colorlinks=true,
	linkcolor=blue,
	hypertexnames=false,
	filecolor=black,      
	urlcolor=blue,
	citecolor=blue
}
\usepackage[nameinlink]{cleveref}

\theoremstyle{plain}
\newtheorem{theorem}{Theorem}[section]
\newtheorem*{theorem*}{Theorem}
\newtheorem{proposition}[theorem]{Proposition}

\theoremstyle{definition}
\newtheorem{definition}[theorem]{Definition}
\newtheorem{construction}[theorem]{Construction}

\theoremstyle{remark}
\newtheorem{remark}[theorem]{Remark}
\newtheorem{example}[theorem]{Example}
\newtheorem*{example*}{Example}

\newcommand{\Z}{\mathbb{Z}}

\newcommand{\keyword}[1]{\textit{\textbf{#1}}}

\newcommand{\framed}[1]{\textrm{Fr}(#1)}

\newcommand{\diag}{\textrm{diag}}
\newcommand{\stab}{\textrm{Stab}}

\thanks{Special thanks to Zachary Greenberg for helpful advice and comments. Part of this work was completed during a summer internship at the MPI MIS sponsored by Anna Wienhard. This project has received funding from the European Research Council (ERC) under the European Union’s Horizon 2020 research and innovation programme (grant agreement No 101018839)
.}

\begin{document}

\begin{abstract}
    We give a precise definition of mutation of skew symmetrizable matrices over group rings and relate it to folding and mutation of quivers with symmetries. These matrices can have non-zero diagonal entries and we explain a mutation rule in some of these cases as well. This new rule comes from a notion of a generalized mutation of an entire quiver or sub-quiver.
\end{abstract}

\maketitle

\section{Introduction}

Quivers and their mutations form the backbone of the theory of cluster algebras since their introduction by Fomin and Zelevinsky \cite{FOZEL02}.  The quivers which appear in this context are always ones without loops of two-cycles; the data of such a quiver is equivalent to a skew symmetric matrix with integer entries called the \keyword{exchange matrix}. Mutation of quivers is equivalent to \keyword{matrix mutation} of the exchange matrices. 

Slightly more generally, it is important in the theory of cluster algebras to consider \textit{skew symmetrizable} exchange matrices. These matrices are used to construct cluster algebras associated with the non-simply laced Dynkin diagrams. Skew symmetrizable exchange matrices can be described by weighted quivers, see for example \cite{kaufman_cluster_2025}  and can often be made by \keyword{folding} \cite{felikson_cluster_2011}. Folding is an operation on quivers which collects their nodes into disjoint sets and then mutates all the nodes of a set together. This \textit{set mutation} operation only makes sense when the nodes of a folding set are disconnected. 

There is of course no reason to restrict to \textit{integer} valued exchange matrices. For example real valued exchange matrices and their mutations have been considered in \cite{Felikson_Tumarkin_2023}.
In this note we consider skew symmetric exchange matrices with entries in a non-commutative ring $\mathcal{R}$ with anti involution $\sigma$. Our primary example is the group ring $\Z[G]$ of a finite group with the involution $\sigma$ induced by $\sigma(g)= g^{-1}$. In this setting we really want \keyword{$\sigma$-skew-symmetric} matricies meaning that $B=-\sigma(B^T)$.
We show that these matrices arise naturally when there is an action of the group $G$ on a quiver $Q$ and we fold $Q$ by this action, as given by \Cref{Main_theorem_construction}.

There is a new difficultly which arises when trying to mutate these matrices; the diagonal entries $b_{ii}$ may not be zero; $b_{ii} = a-\sigma(a)$ is a perfectly allowable as an entry in a $\sigma$ skew symmetric matrix over $\mathcal{R}$. In \Cref{Main_Theorem_mutation}  we show that mutations of the matrices given by given by \Cref{Main_theorem_construction} agree with set mutations of the folding of $Q$ by $G$, assuming that the diagonal entries are zero. 

\Cref{sec:cyclic_groups} of the paper attempts to understand the situation when these diagonal entries are non-zero. Our strategy relies on the notion of a \keyword{generalized mutation} defined in \Cref{def:generalized_mutation}. This generalized mutation is closely related to the \keyword{Green to Red} mutation sequence also known as the \keyword{Cluster DT transformation} see for example \cite{Keller_mutation}. We treat explicitly the cases $b_{ii} = \omega-\omega^{-1}$ in \Cref{sec:3-cycle,sec:4-cycle} for $\omega$ a generator of $\Z/k\Z$ for $k=3,4$ and $b_{ii} =2( \omega-\omega^{-1})$ for $k=3$ in \Cref{sec:double-3}.

\section{Basic Definitions}

\begin{definition}
    A \keyword{quiver} $Q$ is a finite oriented graph that has vertices and directed edges, but no 2-cycles or loops.
\end{definition} 
We will often label the nodes of the quiver $Q$ by $\{1,\dots ,n\}$, in this way there is a natural action of the permutation group $S_n$ on the vertices of $Q$.

\begin{definition}\label{def:quiver_iso}
    Two quivers $Q_1,Q_2$ are \keyword{isomorphic} if there is a permutation $\sigma\in S_n$ such that $\sigma(Q_1)=Q_2$.
\end{definition}

\begin{definition}
    A \keyword{mutation} of $Q$ at node $k$ transforms $Q$ into $Q'=\mu_{k}(Q)$ by changing the arrows through the following rule: 
    \begin{enumerate}
        \item For every pair of arrows $i\rightarrow k$ and $k\rightarrow j$, add a new arrow $i \rightarrow j$.
        \item Remove any oriented 2-cycles.
        \item Reverse all arrows incident to node $k$.
    \end{enumerate}
\end{definition}
Mutation is an involution, so $\mu_{k}(\mu_{k}(Q))= Q$.

\begin{definition} \label{mutation_equivalent}
    Two quivers $Q$ and $Q'$ are \keyword{mutation equivalent} if there is a sequence of mutations that transforms $Q$ into a quiver isomorphic to $Q'$. The \keyword{mutation class} of a quiver is the set of quivers mutation equivalent to $Q$.
\end{definition}

There are two families of mutation classes of quivers we will reference later, namely those of surface type \cite{fomin_cluster_2008} and those associated to Grassmannians \cite{scott_grassmannians_2006}.

\begin{definition}\label{exchange_matrix}
    Let $Q$ be a quiver with $n$ mutable vertices labeled $1,...,n$. Then the \keyword{exchange matrix} of $Q$, $B^Q \in M_n(\Z)$, is an $n\times n$ matrix with integer coefficients defined by: 
    \begin{equation} \label{eq:1}
    (b_{ij}) =
        \begin{cases}
            l & \text{if there are $l$ arrows from vertex $i$ to vertex $j$}\\
            -l & \text{if there are $l$ arrows from vertex $j$ to vertex $i$}\\
            0 & \text{otherwise}
        \end{cases}       
    \end{equation}  
\end{definition}

Generally exchange matrices appearing in the theory of cluster algebras might not necessarily be skew symmetric, but can be \emph{skew-symmetrizable}, see for example the discussion in \cite{fomin2024introductionclusteralgebraschapters}. 

\begin{definition}
A square integer matrix, $B=(b_{ij})$ is called \emph{skew-symmetrizable} if there exists a diagonal skew-symmetrizing matrix, $D$, with diagonal entries $d_{i}$ such that $DB$ is skew-symmetric i.e. $d_{i}b_{ij} = -d_{j}b_{ji}$ for all $i$ and $j$.
\end{definition}
%definition from https://arxiv.org/pdf/math/0104151 

Mutation of quivers can equivalently be defined in terms of a mutation rule for exchange matrices:
\begin{definition}\label{def:max_and_min}
    Let $a\in \Z$. We write $[a]_{+} := \max(a,0)$ and $[a]_{-}=\min(a,0)$ for the \keyword{positive part} and \keyword{negative part} of $a$ respectively. 
\end{definition}

We note that the negative part of $a$ is actually a negative number. We find this to be a more useful definition since $[a]_-$ is either equal to $a$ or 0. 

\begin{definition}[Mutation of Exchange Matrices]\label{def:matrix_mutation}
    The mutation $\mu_k(B) = B'$ where 
     \begin{equation} \label{eq:2}
        (b'_{ij}) =
            \begin{cases}
            -b_{ij} & \text{for $i = k$ or $j=k$}\\
            b_{ij} + [b_{ik}]_{+}[b_{kj}]_{+} - [b_{ik}]_{-}[b_{kj}]_{-} & \text{otherwise.}\\
            \end{cases}       
        \end{equation} 
\end{definition} 

\begin{definition}
    We write $a\circ b:= [a]_+[b]_+ -[a]_-[b]_-$.
\end{definition}

This gives a convenient way to write the mutation rule as $b'_{ij} = b_{ij} + b_{ik}\circ b_{kj}$.

\subsection{Exchange matrices over group rings}

Let $R $ be a ring and  $ \sigma:R\to R$ be an anti-involution i.e. $\sigma(ab)=\sigma(b)\sigma(a)$.

\begin{definition}
    A \keyword{$\sigma$ skew-symmetric matrix} over $R$ is a matrix $B\in M_n(R)$ such that $b_{ij} = -\sigma(b_{ji})$. 

    A matrix $B \in M_n(R)$ is \keyword{$\sigma$ skew-symmetrizable} if there is a diagonal matrix $D$ with entries $d_{ii}$ in the center of $\mathcal{R}$ which are fixed by $\sigma$ such that $DB$ is $\sigma$ skew-symmetric.
\end{definition}

Our main example of such a ring will be the group ring $\Z[G]$ of a finite group with involution given by $\sigma(g)=g^{-1}$. We will always assume that the involution $\sigma$ is given this way. For technical reasons will will extending this ring to include $||G||^{-1}$. 

For the group ring $\Z[G][||G||^{-1}]$ there is a natural choice of the positive and negative part of an element. 

\begin{definition}
   For $a = \sum a_g g$, $a_g \in \Z$ we write $[a]_\pm = \sum[a_g]_\pm g $. 
\end{definition}

With this we can attempt to define mutations of $\sigma$ skew-symmetrizable matrices over $\Z[G][||G||^{-1}]$ using the exact same formula as before. Our goal is to show that such matrices and their mutations naturally arise from folding quivers.

We will write $\epsilon,\omega,\zeta$ for generators of the cyclic groups $\Z/n\Z$ for $n=2,3,4$ respectively.

\section{Symmetric Folding}
\keyword{Folding} is a natural operation on quivers which has been used to construct skew-symmetrizable matrices from skew-symmetric ones. To fold a quiver, we  group its nodes into disjoint sets and do mutations by mutating all the nodes of a given group all at once. Now we recall the notions of folding (with slightly changed terminology) from \cite{Kaufman-special}.

\begin{definition}
    A \keyword{folding} of the quiver $Q$ with nodes $1,...,n$ is a choice of $k$ disjoint and non-empty sets of nodes called the \keyword{folding sets}, where the union of such sets contains all the nodes of $Q$.
    A folding set is called \keyword{cycle-free} if the nodes of this set form a disconnected sub-quiver. A \keyword{set mutation} at a cycle free folding set is a sequence of mutations which mutates each node in the set once.
\end{definition}

\begin{figure}
    \centering
    \begin{tikzcd}
    2 \arrow[rd] &   & 3 \arrow[ld] \\
             & 1 &              \\
    4 \arrow[ru] &   & 5 \arrow[lu]
    \end{tikzcd}
     $\xRightarrow[\text{}]{\text{folding}}$
    \begin{tikzcd}[row sep=.2em]
    \color{red}2 \arrow[rd] &   & \color{blue}3 \arrow[ld] \\
             & 1 &              \\
    \color{red}4 \arrow[ru] &   & \color{blue}5 \arrow[lu]
    \end{tikzcd}
    \caption{Folding Q into 3 folding sets, namely \{1\}, \{2,4\} and \{3,5\}.}
\end{figure}
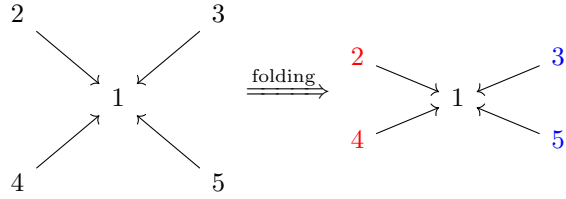

It is natural to define mutation of the folded quiver at cycle free sets by set mutation. Since there are no loops or arrows between the nodes in this, it does not matter which node is mutated first, as long as every node in the set is mutated only once. This is the operation used to relate folding of quivers with mutation of skew symmetrizable exchange matrices in \cite{felikson_cluster_2011}.

\subsection{Folding by the action of a group}

Let $Q$ be a quiver on $n$ nodes and let $G\subset S_n$ be a finite group acting on $Q$ by quiver isomorphisms. 
\begin{definition}
    The folded quiver by $G$ denoted by $Q_G$ is the folding of $Q$ whose groups are given by the orbits of the action of $G$ on $Q$. 
\end{definition}

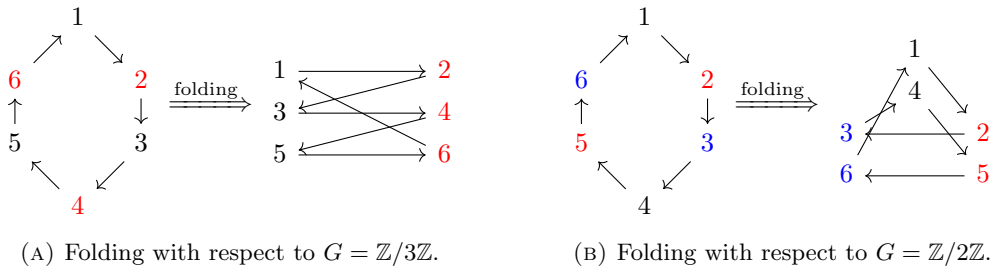
\begin{figure}[b]
    \centering
    \begin{subfigure}{.49\textwidth}
    \centering
    \begin{tikzcd}[sep=1em]
                 & 1 \arrow[rd] &              \\
    \color{red}6 \arrow[ru] &              & \color{red}2 \arrow[d]  \\
    5 \arrow[u]  &              & 3 \arrow[ld] \\
                 & \color{red}4 \arrow[lu] &             
    \end{tikzcd}
    $\xRightarrow{\text{folding}}$
    \begin{tikzcd}[row sep=.2em]
    1 \arrow[rr] &  & \color{red}2 \arrow[lld]  \\
    3 \arrow[rr] &  & \color{red}4 \arrow[lld]  \\
    5 \arrow[rr] &  & \color{red}6 \arrow[lluu]
    \end{tikzcd} 
    \caption{Folding with respect to $G=\Z/3\Z$.}
    %\label{fig:placeholder}
    \end{subfigure}
    \begin{subfigure}{.49\textwidth}
    \centering
    \begin{tikzcd}[sep=1em]
                 & 1 \arrow[rd] &              \\
    \color{blue}6 \arrow[ru] &              & \color{red}2 \arrow[d]  \\
    \color{red}5 \arrow[u]  &              & \color{blue}3 \arrow[ld] \\
                 & 4 \arrow[lu] &             
    \end{tikzcd}
    $\xRightarrow[\text{}]{\text{folding}}$
    \begin{tikzcd}[row sep=.2em,column sep=.6em]                &  & 1 \arrow[rrdd] &  &                           \\
                            &  & 4 \arrow[rrdd] &  &                           \\
\color{blue}3 \arrow[rru]   &  &                &  & \color{red}2 \arrow[llll] \\
\color{blue}6 \arrow[rruuu] &  &                &  & \color{red}5 \arrow[llll]
\end{tikzcd}
        \caption{Folding with respect to $G=\Z/2\Z$.}
    \end{subfigure}
        \caption{Folding the same quiver with respect to two different groups.}
        \label{fig:twofoldings}
\end{figure}

Our goal now is to define a $\sigma$ skew-symmetrizable exchange matrix with entries in the group ring $\Z[G]$ which describes the folded quiver $Q_G$. We denote the orbits of $G$ on $Q$ by the sets $O_1\dots,O_m$

\begin{construction}

Choose a representative for each orbit, $x_{k}\in O_k$ and write $H_k := \stab(x_{k})$ for the stabilizer subgroup of $x_{j}$.
For $g\in G$ we denote by $B_{x_ix_j}^Q(g)$ the number of arrows from node $x_i$ to node $g\cdot x_j$.

The \keyword{folded exchange matrix of $Q$ by $G$}, $B^Q_G(x_1,\dots,x_m) = [c_{ij}]$ is a $m \times m$ matrix with entries given by
    \begin{equation} \label{formula_folded_exchange matrix_entries}
        c_{ij} = \frac{1}{||H_j||}\sum_{g\in G} B^Q_{x_ix_j}(g)g.
    \end{equation}
\end{construction}

\begin{theorem}\label{Main_theorem_construction}
    The matrix $B_G^Q(x_1,\dots,x_m)$ is a skew symmetrizable matrix over the group ring $\Z[G]$.
\end{theorem}
\begin{proof}
     We can see that an arrow from $x_i$ to $gx_j$ implies due to $G$ symmetry that there is an arrow from $g^{-1}x_i$ to $x_j$. This shows that $B^Q_{x_ix_j}(g)=-B^Q_{x_jx_i}(g^{-1})$ and it follows that $$c_{ji}= -\frac{||H_j||}{||H_i||}\sigma(c_{ij}).$$ Thus this matrix is skew symmetrizable with symmetrizer given by $d_{ii}= ||H_i||^{-1}$
\end{proof}

The formula for matrix mutation \Cref{def:matrix_mutation} makes perfect sense for these matrices as well. 

\begin{theorem}\label{Main_Theorem_mutation}
    Let there be a quiver, $Q$, and a group, $G$ acting on $Q$. Suppose that the folding of $Q$ by $G$ satisfies that the orbit $O_j$ is cycle free i.e. the entry $b_{jj}$ of $B_G^Q(x_1,\dots,x_m)$ is zero. Then $\mu_j(B_G^Q(x_1,\dots,x_m))$ agrees with the $B_G^{Q'}$, where $Q'$ is the quiver obtained by set mutation at the orbit $O_j$. 
\end{theorem}

\begin{proof}
    Let $j$ denote the group we mutate at. Clearly the entries of the mutated matrix with one index equal to $j$ change correctly. The other entries $c'_{ik}= c_{ik}+[c_{ij}]_+[c_{jk}]_+-[c_{ij}]_-[c_{jk}]_-$. We write the product of the two positive parts as \begin{align*}
        [c_{ij}]_+[c_{jk}]_+&= \frac{1}{||H_j||\cdot||H_k||}\sum_{h\in G}[B^Q_{x_ix_j}(h)]_+h\sum_{g\in G}[B^Q_{x_jx_k}(g)]_+g \\
        &=\frac{1}{||H_j||\cdot||H_k||}\sum_{h,g\in G}[B^Q_{x_ix_j}(h)]_+[B^Q_{x_jx_k}(h^{-1}g)]_+g \\
        &= \frac{1}{||H_k||}\sum_{hH_j \in G/H_j,g\in G}[B^Q_{x_ix_j}(h)]_+[B^Q_{x_jx_k}(h^{-1}g)]_+g . 
    \end{align*}
    We can write a similar expression for the product of the negative parts. On the other hand the entry $e_{ik}$ of $B^{Q'}_G(x_1,\dots,x_m)$ can be calculated as follows: whenever there is a pair of arrows $(x_i \to hx_j),(hx_j \to gx_k) $ we will add an arrow $(x_i \to gx_k)$, and similarly for arrows in the opposite direction. Therefore we have 
    \begin{align*}
        e_{ik} = c_{ik}+\frac{1}{||H_k||}\sum_{g\in G}\left(\sum_{hH_j\in G/H_j}B^Q_{x_ix_j}(h)\circ B^Q_{x_jx_k}(h^{-1}g)\right)
    \end{align*}
     which is equal to $c_{ik}'$.
\end{proof}

When the entry $b_{jj}\neq 0 $ there is not a yet a well defined notion of set mutation, we will come back to this situation later.

\subsection{Weaving}

The folded exchange matrix $B_G^Q$ depends on the choice of representative of each orbit of $G$. 

\begin{proposition}
    \sloppy The folded exchange matrix $B_G^Q(x_1,\dots,gx_j,\dots,x_m)$ is obtained from $B_G^Q(x_1,\dots,x_m)$ by conjugation by the diagonal matrix $D_j(g)=\diag(1,\dots,g,\dots,1)$ with $g$ in the $j^{th}$ position.
\end{proposition}

\begin{proof}
    $B^Q_{x_igx_j}(h)$ is the number of arrows from $x_i$ to $hgx_j$ which is equal to $B^Q_{x_ix_j}(hg)$ while we also have $B^Q_{gx_jx_i}(h)=B^Q_{x_jx_i}(g^{-1}h)$. Therefore the entry $e_{ij}$ of $B_G^Q(x_1,\dots,gx_j,\dots,x_m)$ is given by 
    $$ e_{ij} = \frac{1}{||H_j||}\sum_{h\in G} B_{x_igx_j}(h)h =  \frac{1}{||H_j||}\sum_{h\in G} B_{x_ix_j}(hg)h = \frac{1}{||H_j||}\sum_{h\in G} B_{x_ix_j}(h)hg^{-1}$$ 
    and the entry $$e_{ji}= \frac{1}{||H_i||}\sum_{h\in G} B_{gx_jx_i}(h)h = \frac{1}{||H_i||}\sum_{h\in G} B_{x_jx_i}(h)gh.$$ 
    This implies that $B_G^Q(x_1,\dots,gx_j,\dots,x_m) = D_j(g)B_G^Q(x_1,\dots,x_j,\dots,x_m)D_j(g)^{-1}$
\end{proof}

\begin{definition}
    We call the change of representative from $x_i$ to $gx_i$ a \keyword{weave} at the group $i$ by the group element $g$. Two matrices $B$ and $B'$ over the group ring $\Z[G]$ are said to be \keyword{weaving isomorphic} if they are conjugate to each other by the product of a diagonal matrix with entries given by group elements $d_{ii}=g_i$ with a permutation matrix.
\end{definition}

We can think of weaving as an action of $G$ on just the nodes in the orbit $i$. Weaving and mutation of quivers commute with each other.

\begin{example}
    Let $G = S_{3}$.  Consider the element $(123)$ of $S_{3}$. The weaving at set $i$, will give us: 

    \begin{equation*}
    \begin{tikzcd}
        i_{1} \arrow[rr]  &  & j_{1} &                                                     &    & i_{1} \arrow[rrdd] &  & j_{1} \\
        i_{2} \arrow[rr] &  & j_{2} & {} \arrow[r, "{Weave}", Rightarrow, maps to] & {} & i_{2} \arrow[rru]  &  & j_{2} \\
        i_{3} \arrow[rr] &  & j_{3} &                                                     &    & i_{3} \arrow[rru]  &  & j_{3}
    \end{tikzcd}   
    \end{equation*}
\end{example}

\subsection{Unfolding}

\begin{definition}
    A \keyword{symmetric unfolding} of a $\sigma$ skew-symmetrizable matrix $B$ over the group ring $\Z[G]$ is a quiver $Q$ with an action of $G$ for which the folded exchange matrix is $B$.
\end{definition}

If we restrict to skew symmetric matrices then we can always construct an unfolding called the \keyword{canonical unfolding}:

\begin{construction}
    Suppose $B$ is an $m\times m$ matrix and let $n$ be the size of $G$. We make a quiver $Q$ with $mn$ nodes labeled $1_{g_1}\dots,1_{g_n},2_{g_1},\dots,2_{g_n},\dots m_{g_n}$. Now for each entry $b_{ij}$ we take the positive part and write it $[b_{ij}]_+ = \sum_k a_{g_k} g_k$ and attach $a_{g_k}$ arrows from node $i_h$ to $j_{hg_k}$ for all $h\in G$. 
\end{construction}
It is easy to see that $G$ acts on this quiver by $g\cdot i_h = i_{gh}$ and the folding by this action taking $1_e, \dots, m_e$ as orbit representatives recovers the matrix $B$. 

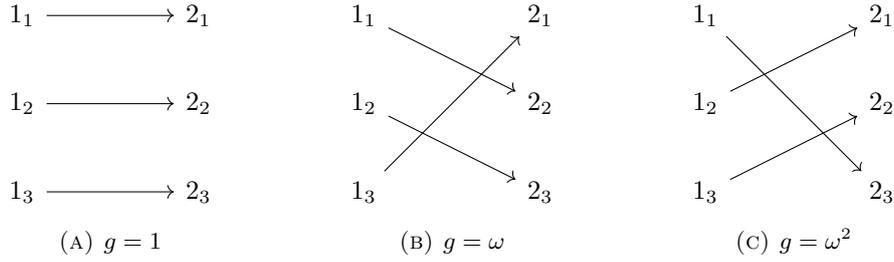
\begin{figure}
    \centering
    \begin{subfigure}{.3\textwidth}
        \centering
        \begin{tikzcd}
        1_{1} \arrow[rr] &    & 2_{1} \\
        1_{2} \arrow[rr] &    & 2_{2} \\
        1_{3} \arrow[rr] &    & 2_{3}
        \end{tikzcd}
        \caption{$g=1$}
    \end{subfigure}
    \begin{subfigure}{.3\textwidth}
        \centering
        \begin{tikzcd}
        1_{1} \arrow[rrd]  &          & 2_{1} \\
        1_{2} \arrow[rrd]  &          & 2_{2} \\
        1_{3} \arrow[rruu] &          & 2_{3}
        \end{tikzcd}
        \caption{$g= \omega$}
    \end{subfigure}
    \begin{subfigure}{.3\textwidth}
        \centering
        \begin{tikzcd}
        1_{1} \arrow[rrdd] &            & 2_{1} \\
        1_{2} \arrow[rru]  &            & 2_{2} \\
        1_{3} \arrow[rru]  &            & 2_{3}
        \end{tikzcd}
        \caption{$g= \omega^2$}
    \end{subfigure}
    \caption{Combinations of arrows of Q with two folding sets of 3 nodes with assigned elements of $G=\Z/3\Z = \{ 1, \omega, \omega^2 \}$.}
    %\label{fig:placeholder}
\end{figure}

\subsection{Examples}

\begin{example}
    Consider the quiver $Q$ and associated exchange matrix: 
    \begin{equation*}
    Q:
        \begin{tikzcd}
        1_{1} \arrow[r] & 2_{1} \arrow[rd] & 3_{1} \arrow[l] \\
        1_{2} \arrow[r] & 2_{2} \arrow[ru] & 3_{2} \arrow[l]
        \end{tikzcd}\quad
    B^Q = 
     \begin{blockarray}{ccccccc}
         & 1_{1} & 1_{2} & 2_{1} & 2_{2} & 3_{1} & 3_{2}\\
       \begin{block}{c[cccccc]}
         1_{1} & 0 & 0 & 1 & 0 & 0 & 0\\
         1_{2} & 0 & 0 & 0 & 1 & 0 & 0\\
         2_{1} & -1 & 0 & 0 & 0 & 0 & 1\\
         2_{2} & 0 & -1 & 0 & 0 & 1 & 0\\
         3_{1} & 0 & 0 & 0 & -1 & 0 & 0\\
         3_{2} & 0 & 0 & -1 & 0 & 0 & 0\\
       \end{block}
     \end{blockarray}
    \end{equation*}

    Choosing representatives $1_{1}$, $2_{1}$ and $3_{1}$, we have the following folded exchange matrix:
    \begin{center}
    $B^Q_{\Z/2\Z}(1_{1}, 2_{1},3_{1}) =
        \begin{bmatrix}
        0 & 1 & 0\\
        -1 & 0 & \epsilon-1\\
        0 & 1-\epsilon & 0
        \end{bmatrix}$
    \end{center}

    If instead we choose representatives $1_{2}$, $2_{1}$ and $3_{1}$, we gave the following folded exchange matrix:
    \begin{center}
    $B^Q_{\Z/2\Z}(1_{2}, 2_{1},3_{1}) =
        \begin{bmatrix}
        0 & \epsilon & 0\\
        -\epsilon & 0 & \epsilon-1\\
        0 & 1-\epsilon & 0
        \end{bmatrix}$
    \end{center}
\end{example}

\begin{example}
    The foldings in \Cref{fig:twofoldings} can be represented by the exchange matricies \\
    $\begin{bmatrix}
        0 & 1- \omega^{-1} \\ -1+\omega &0 
    \end{bmatrix}$ and 
    $\begin{bmatrix}
        0 & 1 & -\epsilon \\ -1 & 0& 1\\ \epsilon &-1 & 0  
    \end{bmatrix}$ up to weaving respectively. 
\end{example}

\begin{example}
    Consider the graph with 12 vertices which is the 1-skeleton of the cubeoctahedron. Make this graph into a quiver, $Q$, by orienting each edge so that the square faces are clockwise and the triangular faces are anticlockwise. This quiver has an action of the symmetric group $S_4$ on it since $S_4$ is the symmetry group of the cubeoctahedron. We can calculate $B^Q_G$ for various subgroups $G\subset S_4$. 

    We have 
    \begin{align*}
        B^Q_{<(1234)>} :
    & \begin{bmatrix}
        (1234)-(1432) & 1-(1234) & 0\\
        -1+(1432)& 0 & 1-(1432) \\
        0 & -1+(1234) & (1432)-(1234)
    \end{bmatrix} \\
B^Q_{<(12)(34),(13)(24)>}: 
     &\begin{bmatrix}
        0 & 1-(12)(34) & -1+(13)(24) \\
        -1+(12)(34)&0& 1-(14)(23) \\
        1-(13)(24)&-1+(14)(23) & 0
    \end{bmatrix}\\
B^Q_{<(123)>}: 
    & \begin{bmatrix}
        (123)-(132) & 1 & -(123) & 0 \\
        -1 & 0 & 1+(123) & -1 \\
        (132)& -1-(132) & 0 & 1\\
        0 &1 & -1&  (132)-(123) 
    \end{bmatrix}\\  
B^Q_{<(1234),(13)>}: 
 &\begin{bmatrix}
    (1234)-(1432) & \frac{1}{2}(1-(1234)+(24)-(13)) \\
    -1+(1432) -(24) +(13) & 0 
\end{bmatrix} \\
B^Q_{G}:
    &\begin{bmatrix}
        \frac{1}{2}\left((1234)-(1432)+(1342)-(1243)\right)
    \end{bmatrix}
    \end{align*}

\end{example}

\section{Exchange Matrices and Mutations over Cyclic Groups }\label{sec:cyclic_groups}

Our goal now is to try to understand general mutations of skew symmetrizable matrices over group rings.

If $B$ is such a matrix and $b_{ii}=0$, then \Cref{Main_Theorem_mutation} implies that the matrix mutation rule agrees with the set mutation of any symmetric unfolding of $B$. What should happen if $b_{ii}\neq0$? 
Our plan is to perform a sequence of mutations on this set which should behave like a single mutation and then refold the quiver to see how the matrix $B$ should change. 

This sequence of mutations should be involutive and have the effect of flipping the direction of arrows which are attached to the orbit $i$. A sequence of mutations of a quiver which changes outward pointing arrows to inward pointing ones is called a \textit{green to red sequence.}

\subsection{The green to red mutation sequence}

\begin{definition}
    A \keyword{framed quiver}, $\framed{Q}$, is an extension of a quiver, $Q$, obtained by adding to each vertex $i$, a new frozen node $i'$ and an arrow $i \rightarrow i'$. We call these frozen nodes the \keyword{framing nodes}. The set of quivers obtained by mutating $\framed{Q}$ at non-frozen nodes is called the \keyword{framed mutation class} of $Q$. Two quivers in the framed mutation class of $Q$ are \keyword{framed isomorphic} if there there is a permutation which fixes the frozen nodes realizing a quiver isomorphism between them. 
\end{definition}

\begin{remark}\label{rem:cluster_framing}
    The elements in the framed mutation class are in correspondence with the clusters in the cluster algebra associated to the quiver $Q$. This provides a convenient way to study the cluster combinatorics of the cluster algebra associated with $Q$ without actually computing any cluster variables or exchange relations.
\end{remark}

After performing some mutations at non-frozen nodes of $\framed{Q}$, the framing nodes will be connected in a more complicated way to the mutable nodes. 

Let $R$ be a quiver in the framed mutation class of $Q$.
\begin{definition}
    A mutable vertex $i$ is a \keyword{green vertex} of $R$ if there are no arrows $j' \rightarrow i$. If a mutable node has no arrows $i \rightarrow j'$ then it is called \keyword{red}.
\end{definition}

Initially, all mutable vertices in $\framed{Q}$ are green. 
It turns out that all nodes must be green or red, this is due to the famous sign coherence of $c$-vectors, a fact which we will not use.

\begin{definition}\label{def:generalized_mutation}
    A sequence of mutations $P = \mu_{i_{1}}...\mu_{i_{N}}$ is called a \keyword{reddening sequence} if all vertices of the quiver $P(\framed{Q})$ are \keyword{red} after performing the mutations. If $P$ satisfies that there is a framed isomorphism $P^2(\framed{Q}) \simeq \framed{Q}$ then we call this sequence a \keyword{generalized mutation} of the quiver $Q$. 
\end{definition}

\subsection{Cyclic Groups}

We now construct some examples of generalized mutations when $G\simeq\Z/n\Z$ 
We write $\omega$ for the generator of $G$. The simplest skew symmetric element of $\Z[G]$ is $\omega-\omega^{-1}$. When $b_{ii}=k(\omega-\omega^{-1})$ we have that the nodes in the $i^{th}$ set of the canonical unfolding of $B$ are connected in an oriented cycle with $k$ arrows.

The proof of the following is a simple calculation using the surface type cluster algebras of \cite{fomin_cluster_2008} including the notion of tagged triangulations of a punctured surface. We do not review these concepts here as they are only needed for this proposition. 

\begin{proposition}
    Let $Q$ be an oriented cycle of length $n$ consisting of \textit{single} arrows. Then $P = \mu_1\mu_2\dots,\mu_n\mu_{n-2}\mu_{n-3}\dots\mu_1$ along with the quiver isomorphism given by the permutation $(n-1,n)$ is a generalized mutation of $Q$.
\end{proposition}
\begin{proof}
    This quiver is mutation equivalent to an orientation of a $D_n$ Dynkin diagram for which its mutation combinatorics and the green to red element is well understood \cite{GoncharovShenQuantumGeometryModuli2022}. The quivers in the mutation class of type $D_n$ are in correspondence with triangulations of a punctured $n$-gon; the oriented cycle corresponds to the wheel triangulation. The mutation sequence given has the effect of moving to the wheel triangulation with all arcs tagged at the puncture. The green to red transformation is obtained by swapping the tagging at the puncture and rotating the boundary one step clockwise; since the wheel triangulation is invariant under rotations we just need to tag the arcs of the triangulation which is exactly what is achieved by the given mutation sequence. 
\end{proof}

\begin{figure}
\begin{tikzcd}[sep=1em]
1_{1} \arrow[rd]  &               \\
               & 1_{2} \arrow[d]  \\
               & 1_{3} \arrow[ld] \\
1_{4} \arrow[uuu] &              
\end{tikzcd}
$\xRightarrow[\text{at $1_{1}$}]{\text{mutation}}$
\begin{tikzcd}[sep=1em]
1_{1} \arrow[ddd] &                         \\
               & 1_{2} \arrow[d] \arrow[lu] \\
               & 1_{3} \arrow[ld]           \\
1_{4} \arrow[ruu] &                        
\end{tikzcd}
$\xRightarrow[\text{at $1_{2}$}]{\text{mutation}}$
\begin{tikzcd}[sep=1em]
1_{1} \arrow[rd] &                \\
              & 1_{2} \arrow[ldd] \\
              & 1_{3} \arrow[u]   \\
1_{4}            &               
\end{tikzcd}
$\xRightarrow[\text{at $1_{3},1_4$}]{\text{mutation}}$
\begin{tikzcd}[sep=1em]
1_{1} \arrow[rd]  &              \\
               & 1_{2} \arrow[d] \\
               & 1_{3}           \\
1_{4} \arrow[ruu] &             
\end{tikzcd}\\
$\xRightarrow[\text{at $1_{2}$}]{\text{mutation}}$
\begin{tikzcd}[sep=1em]
1_{1} \arrow[rdd] &                           \\
               & 1_{2} \arrow[ldd] \arrow[lu] \\
               & 1_{3} \arrow[u]              \\
1_{4} \arrow[ru]  &                          
\end{tikzcd}
$\xRightarrow[\text{at $1_{1}$}]{\text{mutation}}$
\begin{tikzcd}[sep=1em]
1_{1} \arrow[rd] &                \\
              & 1_{2} \arrow[ldd] \\
              & 1_{3} \arrow[luu] \\
1_{4} \arrow[ru] &               
\end{tikzcd}
$\xRightarrow[\text{$(1_{3}, 1_{4})$}]{\text{permutation}}$
\begin{tikzcd}[sep=1em]
1_{1} \arrow[rd]  &               \\
               & 1_{2} \arrow[d]  \\
               & 1_{4} \arrow[ld] \\
1_{3} \arrow[uuu] &              
\end{tikzcd}
\caption{4-cycle returns to itself after performing the generalized mutation.}
\label{fig:4-cycle_mutation}
\end{figure}
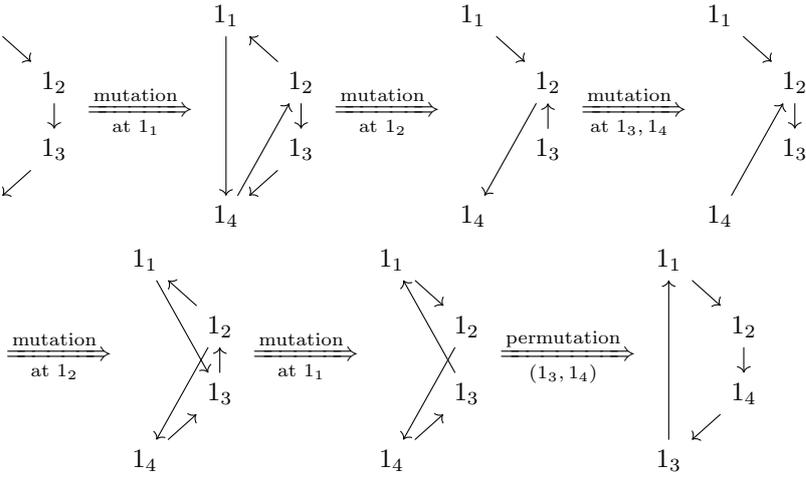

\Cref{fig:4-cycle_mutation} shows the generalized mutation sequence for a 4-cycle.

\begin{remark}
    The mutation sequence for $n=3$ has appeared in several different contexts in the literature, e.g. \cite{baziermatte2024knottheoryclusteralgebras} in a mutation sequence representing the Reidemeister 3 move and in \cite{martello2024okamotossymmetryrepresentationspace} as the mutation realizing a birational transformation related to the Painlevé VI equation. 
\end{remark}

\subsubsection{The 3-Cycle}\label{sec:3-cycle}
Now we use this generalized mutation sequence for $G\simeq\Z/3\Z$. Consider the following exchange matrix: $$    B =
\begin{bmatrix}
        \omega - \omega^{-1} & a + b\omega +c\omega^{-1} & -(d+e\omega^{-1} +f\omega) \\
        -(a + b\omega^{-1} +c\omega) & 0 & 0 \\
        (d+e\omega +f\omega^{-1}) & 0 & 0 
\end{bmatrix}
$$
The canonical unfolding of this matrix is given by the quiver
    \begin{equation*}\label{example_3_cycle_general}
    \begin{tikzcd}
3_{1} \arrow[rr, "d",red, no head] \arrow[rrrd,  no head] \arrow[rrdd,  no head] &  & 1_{1} \arrow[rd] \arrow[rrr, "a",red, no head] \arrow[rrrd, "b",blue, no head] \arrow[rrrdd, "c",orange,no head] &                                                                                  &  & 2_{1} \\
3_{2} \arrow[rrr, no head] \arrow[rru,"e",blue, no head] \arrow[rrd, no head]                 &  &                                                                                           & 1_{2} \arrow[ld] \arrow[rr, no head] \arrow[rru,  no head] \arrow[rrd, no head] &  & 2_{2} \\
3_{3} \arrow[rr, no head] \arrow[rrru, no head] \arrow[rruu,"f",orange, no head]                &  & 1_{3} \arrow[uu] \arrow[rrr, no head] \arrow[rrruu, no head] \arrow[rrru, no head]           &                                                                                  &  & 2_{3}
\end{tikzcd}
\end{equation*}
where $a,b,c$ are integers representing the number of arrows from one node to another, which may be positive or negative.

Our generalized mutation for the set 1 is the sequence of mutations $\mu_{1_1}\mu_{1_2}\mu_{1_3}\mu_{1_1}$ followed by the  permutation exchanging nodes $1_{2}$ and $1_{3}$.

\begin{example}\label{example_mutation_3_cycle}
For example, mutating the quiver  below at set 1 using the rule above has the following result:
    \begin{center}
    \begin{tikzcd}[row sep=1em]
    1_{1} \arrow[rd] \arrow[rr,red] &                         & 2_{1} \\
                         & 1_{2} \arrow[ld] \arrow[r,red] & 2_{2} \\
    1_{3} \arrow[uu] \arrow[rr,red] &                         & 2_{3}
    \end{tikzcd}
    $\xRightarrow[\text{at 1}]{\text{mutation}}$ 
    \begin{tikzcd}[row sep=1em]
    1_{1} \arrow[rd] &                        &  & 2_{1} \arrow[lld, orange]   \\
              & 1_{3} \arrow[ld]          &  & 2_{2} \arrow[llld, orange]  \\
    1_{2} \arrow[uu] &                        &  & 2_{3} \arrow[llluu, orange]
    \end{tikzcd}   
    \end{center}

\end{example}

Calculating the effect of this sequence of mutations for general values of $b_{12},b_{31}$ is simple but tedious. We can collect the calculations together by writing them in terms of group ring elements. Essentially, we will only unfold the set 1 while leaving nodes 2 and 3 unfolded:  
\begin{equation*}
    \begin{tikzcd}
                                                                                                                      &  &  & 1_{1} \arrow[rdd] \arrow[rrrrdd, "b_{12}", no head]              &                                                         &  &  &   \\
                                                                                                                      &  &  &                                                                  &                                                         &  &  &   \\
3 \arrow[rrrr, "\frac{1}{3}\omega^{-1}b_{31}", no head] \arrow[rrruu, "\frac{1}{3}b_{31}", no head] \arrow[rrrdd, "\frac{1}{3}\omega b_{31}"', no head] &  &  &                                                                  & 1_{2} \arrow[ldd] \arrow[rrr, "\omega b_{12}", no head] &  &  & 2 \\
                                                                                                                      &  &  &                                                                  &                                                         &  &  &   \\
                                                                                                                      &  &  & 1_{3} \arrow[uuuu] \arrow[rrrruu, "\omega^{-1} b_{12}"', no head] &                                                         &  &  &  
\end{tikzcd}
\end{equation*}

Now calculating the mutation sequence gives the following result: 

\begin{proposition}
    Suppose $b_{kk}=\omega-\omega^{-1}$. Mutation at $k$ changes $B$ by the following: 
    \begin{equation}\label{eq:cyc_3_aligned}
        (b'_{ij}) =
        \begin{cases}
        b_{ij}
        - [\omega^{-1}b_{ij} + [b_{ij}]_{+}]_{+}
        - [\omega b_{ij} + [b_{ij}]_{-}]_{-},
        & \text{for $i = k$ or $j = k$,} \\[1ex]

        \begin{aligned}
        b_{ij}
        &+ \frac{1}{3}\Big(
        b_{ik}\circ b_{kj}
        - b'_{ik}\circ b'_{kj} \\
        &\qquad
        + (\omega^{-1}b_{ik} + [b_{ik}]_{+})
        \circ (\omega b_{kj} + [b_{kj}]_{-}) \\
        &\qquad
        + (\omega b_{ik} + [b_{ik}]_{-})
        \circ (\omega^{-1}b_{kj} + [b_{kj}]_{+})
        \Big),
        \end{aligned}
        & \text{otherwise.}
        \end{cases}
    \end{equation}
\end{proposition}

If $b_{12}=a+b\omega+c\omega^{-1}$ then writing $b'_{12}=a'+b'\omega+c'\omega^{-1}$
we can calculate that $a' = a - [b+[a]_{+}]_{+} - [c+[a]_{-}]_{-}$. We can analyze each case to find: 
\begin{equation*} \label{eq: 3-node_a'}
    a' =
    \begin{cases}
            \text{if $a>0$}
            \begin{cases}
                \text{if $b+a>0$}
                \begin{cases}
                    \text{if $c>0$} & $=-b$
                    \\
                    \text{if $c<0$} & $=-b-c$
                \end{cases}
                \\
                \text{if $b+a<0$}
                \begin{cases}
                    \text{if $c>0$} & $=a$
                    \\
                    \text{if $c<0$} & $=a-c$
                \end{cases}
            \end{cases}
            \\
            \text{if $a<0$}
            \begin{cases}
                \text{if $c+a>0$}
                \begin{cases}
                    \text{if $b>0$} & $=a-b$
                    \\
                    \text{if $b<0$} & $=a$
                \end{cases}
                \\
                \text{if $c+a<0$}
                \begin{cases}
                    \text{if $b>0$} & $=-b-c$
                    \\
                    \text{if $b<0$} & $=-c$
                \end{cases}
            \end{cases}
    \end{cases}  
\end{equation*}
with similar formulas for $b'$ and $c'$.

\begin{example}
The quiver $Q_{3n}$ is of mutation type $A_{3n}$, corresponding to a triangulation of a $3n+3$-gon. 
\begin{equation*}
        Q_{3n} := \begin{tikzcd}
1_1 \arrow[rd] &                        & 2_1 \arrow[ll] & 3_1 \arrow[l] & \dots \arrow[l] & n_1 \arrow[l] \\
               & 1_2 \arrow[ld]         & 2_2 \arrow[l]  & 3_2 \arrow[l] & \dots \arrow[l] & n_2 \arrow[l] \\
1_3 \arrow[uu] &                        & 2_3 \arrow[ll] & 3_3 \arrow[l] & \dots \arrow[l] & n_3 \arrow[l]
\end{tikzcd}
\end{equation*}
Its symmetric folding by $\Z/3\Z$ is also of finite mutation type, it is easy to see that its mutations correspond to 3-fold symmetric triangulations of a $3n+3$-gon. 
\end{example}

\begin{example}
    The exchange matrices  
    \begin{align*}
         &W_3= \begin{bmatrix}
        0 & 1-\omega & 0 \\
        -1+\omega^{-1} & 0 & 1\\
        0 & -1 &0
    \end{bmatrix}\quad W_4= \begin{bmatrix}
        0 & 1-\omega & 0 &-1\\
        -1+\omega^{-1} & 0 & 1 &0\\
        0 & -1 &0 & 1 \\
        1 & 0 &-1 &0 
    \end{bmatrix}
    \\ &W_5= \begin{bmatrix}
        0 & 1-\omega & 0 &0 &0\\
        -1+\omega^{-1} & 0 & 1 &0 & -1\\
        0 & -1 &0 & -1 & 1+\omega \\
        0 & 0 &1 &0 & -\omega \\
        0 & 1 & -1 -\omega^{-1} & \omega^{-1} & 0 
    \end{bmatrix}  
    \end{align*}
    are all finite mutation type. Using \Cref{rem:cluster_framing} we can compute the exchange graph of a cluster algebra associated to these matrices by framing their canonical unfoldings, as shown in \Cref{fig:exchange_complexes}. The canonical unfoldings of $W_3$ and $W_4$ are in the mutation classes of the Grassmannian cluster algebras $Gr(4,8)$ and $Gr(4,9)$ while the unfolding of $W_5$ is of the mutation class of the triangular quiver associated with the base affine space $\mathbb{C}[SL_7(\mathbb{C})]$ discussed in e.g. \cite{fomin2021introductionclusteralgebraschapter6}. The folding of the Grassmannian cluster algebra $Gr(4,9)$ by $\Z/3\Z$ giving $W_4$ is given by the cyclic symmetry shifting the 9 vectors by $v_i \to v_{i+3}$, see \cite{fraser_cyclic_2020}.

    \begin{figure}
        \centering
        \begin{subfigure}{.45\textwidth}
        \centering
            \includegraphics[scale=.09]{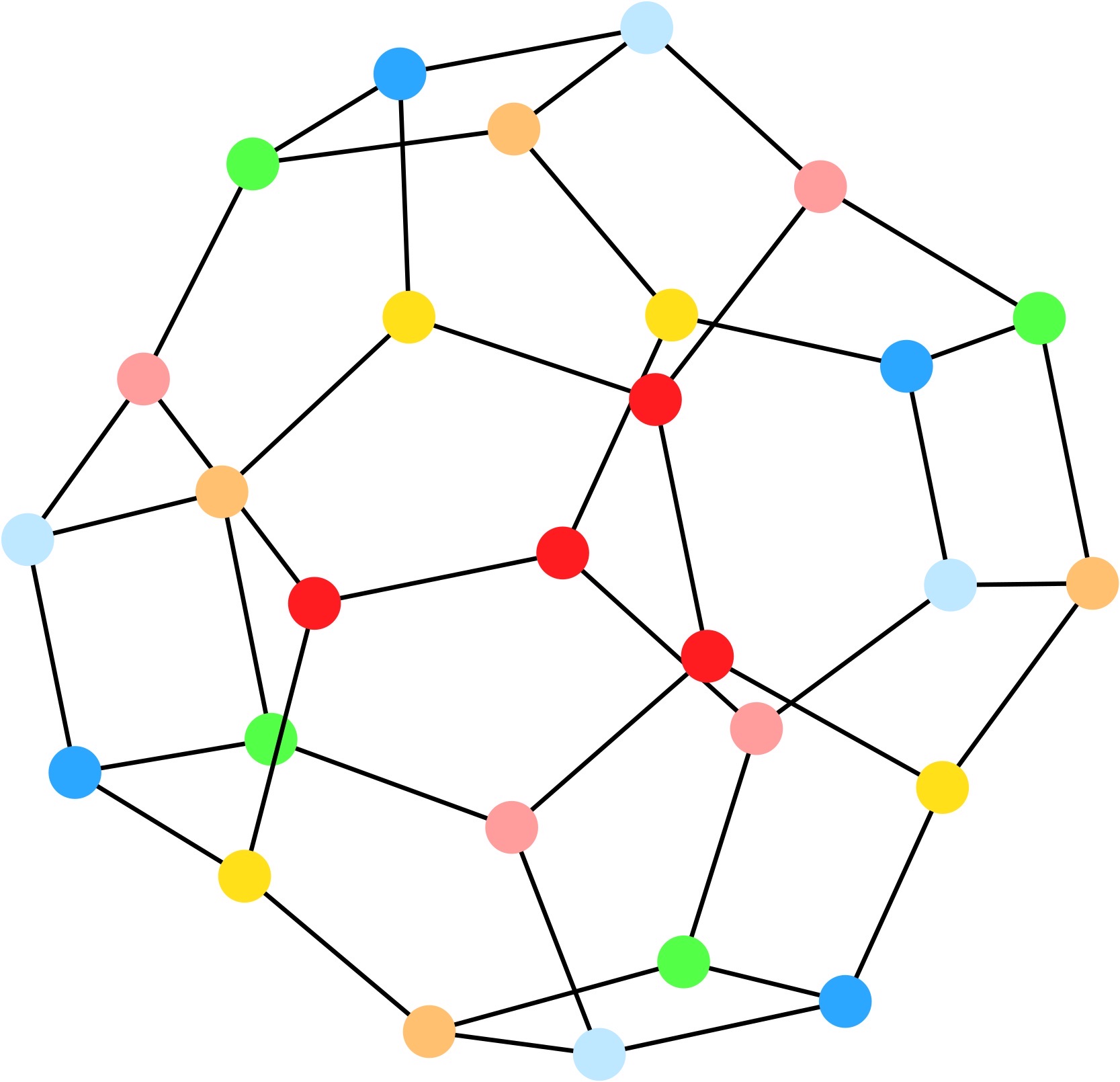}
            \caption{Exchange Graph of $W_3$}
        \end{subfigure}
        \begin{subfigure}{.45\textwidth}
        \centering
            \includegraphics[scale=.05]{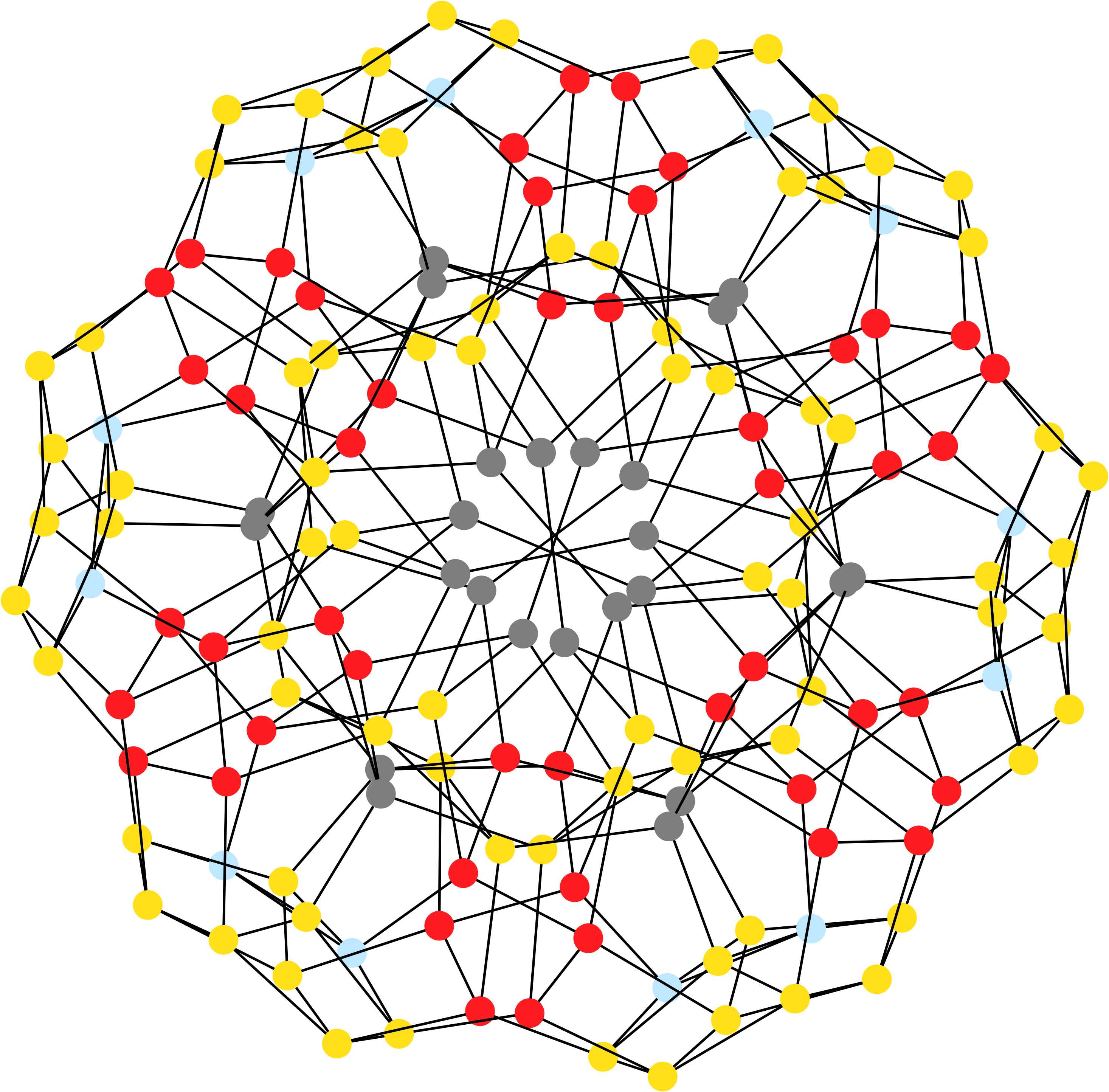}
            \caption{Exchange Graph of $W_4
            $}
        \end{subfigure}
        \label{fig:exchange_complexes}
    \end{figure}
 
\end{example}
%Define cluster variables?

\subsubsection{The 4-Cycle}\label{sec:4-cycle}
The next obvious step is to ask what happens when instead of 3-node cycle we have a cycle of 4 nodes and what the mutation formula is in this case. We consider the exchange matrix $B$:
\begin{equation*}
   B=
\begin{bmatrix}
        \zeta - \zeta^{-1} & a + b\zeta +c\zeta^{2} + d\zeta^{3} & -(e+f\zeta^{3}+g\zeta^{2}+h\zeta)\\
        -(a + b\zeta^{3} +c\zeta^{2} + d\zeta) & 0 & 0  \\
        e+f\zeta+g\zeta^{2}+h\zeta^4 &0&0
\end{bmatrix}
\end{equation*}

where the combinations of arrows can be represented with elements of the group $G = \Z/4\Z = \{1,\zeta, \zeta^{2}, \zeta^{3}\}$. 

Then the canonical unfolding of this matrix is given by the quiver: 
\begin{equation*}
\begin{tikzcd}[column sep=4em,row sep=1em]
3_{1} \arrow[rrr, "e", red, no head] \arrow[rrrrd,blue,"f" no head] \arrow[rrrrdd,"g", orange, no head] \arrow[rrrddd,"h", Green, no head] &  &  & 1_{1} \arrow[rrr, "a", red, no head] \arrow[rd] \arrow[rrrd, "b", blue, no head] \arrow[rrrdd, "c", orange, no head] \arrow[rrrddd, "d", Green, no head] &                                                                                                     &  & 2_{1} \\
3_{2} \arrow[rrrr, no head] \arrow[rrru, no head] \arrow[rrrrd, no head] \arrow[rrrdd, no head]   &  &  &                                                                                                                                & 1_{2} \arrow[rr, no head] \arrow[d] \arrow[rru, no head] \arrow[rrd, no head] \arrow[rrdd, no head] &  & 2_{2} \\
3_{3} \arrow[rrrr, no head] \arrow[rrruu,  no head] \arrow[rrrru, no head] \arrow[rrrd, no head]   &  &  &                                                                                                                                & 1_{3} \arrow[rr, no head] \arrow[ld] \arrow[rruu, no head] \arrow[rru, no head]                     &  & 2_{3} \\
3_{4} \arrow[rrr, no head] \arrow[rrruuu,  no head] \arrow[rrrruu, no head] \arrow[rrrru, no head] &  &  & 1_{4} \arrow[rrr, no head] \arrow[uuu] \arrow[rrru, no head] \arrow[rrruu, no head] \arrow[rrruuu, no head]                    &                                                                                                     &  & 2_{4}
\end{tikzcd}
\end{equation*}

Similarly to the case of 3 nodes in set 1, we need a sequence of mutations and permutation that will return the cycle back to itself, which will correspond to mutating at the set of nodes 1. The sequence is in this case: $\mu_{1}\mu_{2}\mu_{4}\mu_{3}\mu_{2}\mu_{1}$  and the permutation exchanges $1_{3}$ and $1_{4}$.

\begin{example}\label{example_mutation_4_cycle}
For example, mutating the quiver below at set 1 using the sequence above has the following effect: 
\begin{center}
    \begin{tikzcd}[row sep=.5em]
    1_1 \arrow[rrrr, red] \arrow[rd]  &                           &                   &  & 2_{1} \\
                            & 1_{2} \arrow[rrr, red] \arrow[d]  &                   &  & 2_{2} \\
                            & 1_{3} \arrow[rrr, red] \arrow[ld] &                   &  & 2_{3} \\
    1_{4} \arrow[rrrr, red] \arrow[uuu] &                           &                   &  & 2_{4}
    \end{tikzcd}
    $\xRightarrow[\text{at 1}]{\text{mutation}}$
    \begin{tikzcd}[row sep=.5em]
    1_1 \arrow[rd]  &               &                   &  & 2_{1} \arrow[llld,cyan]    \\
               & 1_{2} \arrow[d]  &                   &  & 2_{2} \arrow[llld, cyan]    \\
               & 1_{3} \arrow[ld] &                   &  & 2_{3} \arrow[lllld, cyan]   \\
    1_{4} \arrow[uuu] &               &                   &  & 2_{4} \arrow[lllluuu,bend right,cyan]
    \end{tikzcd}
    
\end{center}
\end{example}
Folding nodes 2 and 3 and leaving 1 unfolded, we have: 
\begin{equation*}
\begin{tikzcd}
                                                                                                                                                                                                          &  &  & 1_1 \arrow[rrrrrdd, "b_{12}", no head] \arrow[rd]              &                                                           &  &  &  &   \\
                                                                                                                                                                                                          &  &  &                                                                & 1_{2} \arrow[rrrrd, "\zeta b_{12}"', no head] \arrow[dd]  &  &  &  &   \\
3 \arrow[rrruu, "\frac{1}{4}b_{31}", no head] \arrow[rrrru, "\frac{1}{4}\zeta^{3}b_{31}", no head] \arrow[rrrrd, "\frac{1}{4}\zeta^{2}b_{31}"', no head] \arrow[rrrdd, "\frac{1}{4}\zeta b_{31}"', no head] &  &  &                                                                &                                                           &  &  &  & 2 \\
                                                                                                                                                                                                          &  &  &                                                                & 1_{3} \arrow[rrrru, "\zeta^2 b_{12}", no head] \arrow[ld] &  &  &  &   \\
                                                                                                                                                                                                          &  &  & 1_{4} \arrow[rrrrruu, "\zeta^3 b_{12}"', no head] \arrow[uuuu] &                                                           &  &  &  &  
\end{tikzcd}
\end{equation*}

\begin{proposition}\label{eq:cyc_4_aligned}
    Suppose $b_{kk}=\zeta-\zeta^{-1}$. Then the mutation at $k$ changes $B$ by the following: 
    \begin{comment}
    \begin{equation}\label{eq:cyc_4_aligned}
        (b'_{ij}) =
        \begin{cases}
        b_{ij}
        - [b_{ij}\zeta + [b_{ij}]_{-}]_{-} - [-b_{ij}\zeta - [b_{ij}]_{-} + [b_{ij}\zeta^{2} + [b_{ij}\zeta + [b_{ij}]_{-}]_{-}]_{-} + [b_{ij}\zeta^3 + [b_{ij}]_{+} + [b_{ij}\zeta+[b_{ij}]_{-}]_{+}]_{+}]_{+},
        & \text{for $i = k$ or $j = k$,} \\[1ex]

        b_{ij} + \frac{1}{4} [b_{ik} \circ b_{kj} - b'_{ik} \circ b'_{kj} 
        + (b_{ik}\zeta^{-1} + [b_{ik}]_{+})\circ(b_{kj}\zeta + [b_{kj}]_{-}) 
        + (b_{ik}\zeta^2 + [b_{ik}\zeta^{-1} + [b_{ik}]_{+}]_{+})\circ(b_{kj}\zeta^{2} + [b_{kj}\zeta + [b_{kj}]_{-}]_{-}) 
        + (b_{ik}\zeta + [b_{ik}]_{-} + [b_{ik}\zeta^{-1} + [b_{ik}]_{+}]_{-})\circ(b_{kj}\zeta^{-1} + [b_{kj}]_{+} + [b_{kj}\zeta+ [b_{kj}]_{-}]_{+}) 
        + (-b_{ik}\zeta^{-1} - [b_{ik}]_{+} + [b_{ik}\zeta^{2} + [b_{ik}\zeta^{-1} + [b_{ik}]_{+}]_{+}]_{+} + [b_{ik}\zeta + [b_{ik}]_{-} + [b_{ik}\zeta^{-1} + [b_{ik}]_{+}]_{-}]_{-})\circ(-b_{kj}\zeta - [b_{kj}]_{-} + [b_{kj}\zeta^{2} + [b_{kj}\zeta + [b_{kj}]_{-}]_{-}]_{-} + [b_{kj}\zeta^{-1} + [b_{kj}]_{+} + [b_{kj}\zeta + [b_{kj}]_{-}]_{+}]_{+})]& \text{otherwise.}
        \end{cases}
    \end{equation}
    \end{comment}
    For $i=k$ or $j=k$:
    \begin{equation*}b_{ij}'=
        \begin{aligned}
        b_{ij}
        &- [b_{ij}\zeta + [b_{ij}]_{-}]_{-} \\
        &- \bigl[
            -b_{ij}\zeta - [b_{ij}]_{-}
            + [b_{ij}\zeta^{2} + [b_{ij}\zeta + [b_{ij}]_{-}]_{-}]_{-} \\
        &
        + [b_{ij}\zeta^3 + [b_{ij}]_{+}
         + [b_{ij}\zeta+[b_{ij}]_{-}]_{+}]_{+}
        \bigr]_{+},
        \end{aligned}
        \end{equation*}
        Otherwise we have: $b_{ij}'=$
        \begin{equation*}
        \begin{aligned}
        b_{ij}
        &+ \frac{1}{4}\Bigl[
         b_{ik} \circ b_{kj}
        - b'_{ik} \circ b'_{kj} \\
        &
            + (b_{ik}\zeta^{-1} + [b_{ik}]_{+})
            \circ (b_{kj}\zeta + [b_{kj}]_{-}) \\
        &
            + (b_{ik}\zeta^2 + [b_{ik}\zeta^{-1} + [b_{ik}]_{+}]_{+})
            \circ (b_{kj}\zeta^{2} + [b_{kj}\zeta + [b_{kj}]_{-}]_{-}) \\
        &
         + (b_{ik}\zeta + [b_{ik}]_{-} + [b_{ik}\zeta^{-1} + [b_{ik}]_{+}]_{-})\circ(b_{kj}\zeta^{-1} + [b_{kj}]_{+} + [b_{kj}\zeta+ [b_{kj}]_{-}]_{+}) \\
        &
         + (-b_{ik}\zeta^{-1} - [b_{ik}]_{+} + [b_{ik}\zeta^{2} + [b_{ik}\zeta^{-1} + [b_{ik}]_{+}]_{+}]_{+} + [b_{ik}\zeta + [b_{ik}]_{-} + [b_{ik}\zeta^{-1} + [b_{ik}]_{+}]_{-}]_{-}) \\
        &
         \circ(-b_{kj}\zeta - [b_{kj}]_{-} + [b_{kj}\zeta^{2} + [b_{kj}\zeta + [b_{kj}]_{-}]_{-}]_{-} + [b_{kj}\zeta^{-1} + [b_{kj}]_{+} + [b_{kj}\zeta + [b_{kj}]_{-}]_{+}]_{+})
        \Bigr],
        \end{aligned}
    \end{equation*}
\end{proposition} 

If $b_{12} =a+b\zeta+c\zeta^{2}+d\zeta^{3}$ then $b_{12}' = a'+b'\zeta+c'\zeta^{2}+d'\zeta^{3}$ we find 
$a' = a - [d + [a]_{-}]_{-} - [-d -[a]_{-} + [c+[d+[a]_{-}]_{-}]_{-} + [b + [a]_{+} + [d+[a]_{-}]_{+}]_{+}]_{+}$ and we can analyze each case to get: 
\begin{equation*}\label{eq: 4-node_a'}
    a' =
    \begin{cases}
            \text{if $a>0$}
            \begin{cases}
                \text{if $d>0$}
                \begin{cases}
                    \text{if $c>0$}
                    \begin{cases}
                        \text{if $b+a+d>0$} & =-b
                        \\
                        \text{if $b+a+d<0$} & =a
                    \end{cases}
                    \\
                    \text{if $c<0$}
                    \begin{cases}
                        \text{if $b+a+d>0$} & $$=a - (b+a+c)_{+}$$
                        \\
                        \text{if $b+a+d<0$} & =a
                    \end{cases}
                \end{cases}
                \\
                \text{if $d<0$}
                \begin{cases}
                    \text{if $c+d>0$}
                    \begin{cases}
                        \text{if $b+a>0$} & =-b
                        \\
                        \text{if $b+a<0$} & =a
                    \end{cases}
                    \\
                    \text{if $c+d<0$}
                    \begin{cases}
                        \text{if $b+a>0$} & $$=a-d-(c+a+b)_{+}$$
                        \\
                        \text{if $b+a<0$} & =a-d
                    \end{cases}
                \end{cases}
            \end{cases}
            \\
            \text{if $a<0$}
            \begin{cases}
                \text{if $a+d>0$}
                \begin{cases}
                    \text{if $b+a+d>0$}
                    \begin{cases}
                        \text{if $c>0$} & =a-b
                        \\
                        \text{if $c<0$} & =a-(c+b)_{+}
                    \end{cases}
                    \\
                    \text{if $b+a+d<0$}
                    \begin{cases}
                        \text{if $c>0$} & =a
                        \\
                        \text{if $c<0$} & =a
                    \end{cases}
                \end{cases}
                \\
                \text{if $a+d<0$}
                \begin{cases}
                    \text{if $b>0$}
                    \begin{cases}
                       \text{if $c+d+a>0$} & =a-b
                       \\
                       \text{if $c+d+a<0$} & =-d-(c+b)_{+}
                    \end{cases}
                    \\
                    \text{if $b<0$}
                    \begin{cases}
                        \text{if $c+d+a>0$} & =a
                        \\
                        \text{if $c+d+a<0$} & =-d
                    \end{cases}
                \end{cases}
            \end{cases}
    \end{cases}  
\end{equation*}

\begin{remark}
    The first formula in \Cref{eq:cyc_4_aligned} appears by considering the arrows between nodes $1_1$ and 2 after performing the mutation sequence. Symmetrically calculating the arrows between the other nodes $1_i$ and 2 will give the same result; however this gives very different formula. It may be possible to simplify the expression using these inequivalent forms. 
\end{remark}

\begin{example}
     The exchange matrices  
    \begin{equation*}
        U_3= \begin{bmatrix}
        0 & 1-\zeta & 0 \\
        -1+\zeta^{-1} & 0 & 1\\
        0 & -1 &0
    \end{bmatrix}\quad U_4= \begin{bmatrix}
        0 & 1-\zeta & 0 &0\\
        -1+\zeta^{-1} & 0 & 1 &-1\\
        0 & -1 &0 & 1+\zeta \\
        0 & -1 &-1 -\zeta^{-1} &0 
    \end{bmatrix}
    \end{equation*}
    are finite mutation type. The latter $U_4$ has only 7 elements in its mutation class up to weaving isomorphism. Its canonical unfolding is mutation equivalent to the mutable portion of the quiver underlying the Grassmannian cluster algebra $Gr(5,10)$.
\end{example}

\subsubsection{The Double 3-Cycle}\label{sec:double-3}
Instead of increasing the number of nodes in the cycle, we can also increase the number of arrows connecting each node in the loop. In this case it is not clear if we can find a generalized mutation. 

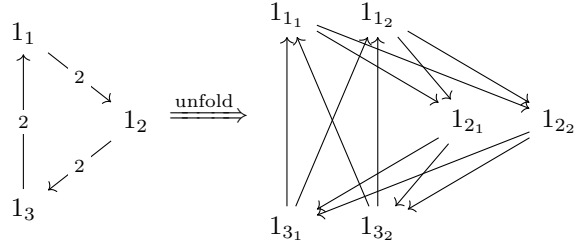
\begin{figure}
    \centering
    \begin{tikzcd}
1_1 \arrow[rd, "2" description] &                                \\
                               & 1_{2} \arrow[ld, "2" description] \\
1_{3} \arrow[uu, "2" description] &                               
\end{tikzcd}
$\xRightarrow[]{\text{unfold}}$
\begin{tikzcd}[sep=1.2em]
1_{1_{1}} \arrow[rrdd] \arrow[rrrdd] & 1_{1_{2}} \arrow[rrdd] \arrow[rdd]   &                                 &                                   \\
                                  &                                   &                                 &                                   \\
                                  &                                   & 1_{2_{1}} \arrow[lldd] \arrow[ldd] & 1_{2_{2}} \arrow[llldd] \arrow[lldd] \\
                                  &                                   &                                 &                                   \\
1_{3_{1}} \arrow[uuuu] \arrow[ruuuu] & 1_{3_{2}} \arrow[uuuu] \arrow[luuuu] &                                 &                                  
\end{tikzcd}

    \caption{Unfolding the double 3-cycle}
    \label{fig:unfolding_3-cycle}
\end{figure}

The double 3-cycle is the well known \textit{Markov Quiver} which is of surface type associated to a once punctured torus, and is there is known to have no green-to-red mutation sequence \cite{Keller_mutation}. In order to solve this problem, we can unfold the 3-cycle, as shown below, and find a generalized mutation for the unfolded quiver which we will then fold back to the original double 3-cycle quiver, see \Cref{fig:unfolding_3-cycle}.

The unfolded 3-cycle is the quiver associated to a 4-punctured sphere, this quiver has a red to green element of order two  which 
is reminiscent to that of the 3-node cycle case with one arrow between each node: $$(\mu_{{1_{1}}_{2}}\mu_{{1_{3}}_{2}}\mu_{{1_{3}}_{1}}\mu_{{1_{2}}_{2}}\mu_{{1_{2}}_{1}}\mu_{{1_{1}}_{2}})(\mu_{{1_{1}}_{1}}\mu_{{1_{3}}_{2}}\mu_{{1_{3}}_{1}}\mu_{{1_{2}}_{2}}\mu_{{1_{2}}_{1}}\mu_{{1_{1}}_{1}})$$

For simplicity we will only calculate the change of adjacent entries of the exchange matrix.
Let $B =\begin{bmatrix}
        2(\omega - \omega^{-1}) & a + b\omega +c\omega^{-1}  \\
        -(a + b\omega^{-1} +c\omega) & 0  \\
\end{bmatrix}$. 
 We unfold the first group further:
\begin{equation*}
    \begin{tikzcd}
1_1 \arrow[rdd, "2" description] \arrow[rrdd, "b_{12}", no head]               &                                                                      &   \\
                                                                               &                                                                      &   \\
                                                                               & 1_2 \arrow[ldd, "2" description] \arrow[r, "\omega b_{12}", no head] & 2 \\
                                                                               &                                                                      &   \\
1_3 \arrow[uuuu, "2" description] \arrow[rruu, "\omega^{-1} b_{12}"', no head] &                                                                      &  
\end{tikzcd}
    \xRightarrow[]{\text{unfold}}
\begin{tikzcd}[row sep=1.1em]
1_{1_{1}} \arrow[rrddd] \arrow[rrdddd] \arrow[rrrrrrrddd, "b_{12}", no head]              &                                                                                             &                                                                               &  &  &  &  &   \\
                                                                                          & 1_{1_{2}} \arrow[rddd] \arrow[rdd] \arrow[rrrrrrdd, "b_{12}"', no head]                      &                                                                               &  &  &  &  &   \\
                                                                                          &                                                                                             &                                                                               &  &  &  &  &   \\
                                                                                          &                                                                                             & 1_{2_{1}} \arrow[llddd] \arrow[ldddd] \arrow[rrrrr, "\omega b_{12}", no head] &  &  &  &  & 2 \\
                                                                                          &                                                                                             & 1_{2_{2}} \arrow[lldd] \arrow[lddd] \arrow[rrrrru, "\omega b_{12}", no head]  &  &  &  &  &   \\
                                                                                          &                                                                                             &                                                                               &  &  &  &  &   \\
1_{3_{1}} \arrow[uuuuuu] \arrow[ruuuuu] \arrow[rrrrrrruuu, "\omega^{-1} b_{12}"', no head] &                                                                                             &                                                                               &  &  &  &  &   \\
                                                                                          & 1_{3_{2}} \arrow[uuuuuu] \arrow[luuuuuuu] \arrow[rrrrrruuuu, "\omega^{-1} b_{12}"', no head] &                                                                               &  &  &  &  &  
\end{tikzcd}
\end{equation*}
Now we perform the mutation sequence and calculate $b_{12}$. The intemidiate quantities involved in the calculation are:
\begin{align*}
    &C:= \omega b_{12}+[b_{12}]_- &D:= \omega^{-1}b_{12}+[b_{12}]_+  \\
    &E:= -b_{12}+2[C]_-+2[D]_+ &F:= b_{12}+ 2[C]_++2[D]_-  \\
    &G:= -C+[E]_-+[F]_+ &H:= -D+[E]_++[F]_- 
\end{align*}
Finally we find $b_{12}' = -E + 2[G]_-+2[H]_+$.

\printbibliography

\end{document}